\documentclass[runningheads,envcountsame,a4paper]{llncs} 

\usepackage{amssymb,amsmath}
\usepackage{wasysym}
\usepackage{graphicx}
\usepackage{epsfig}
\usepackage{latexsym}
\usepackage[english]{babel}
\usepackage[utf8]{inputenc}
\usepackage{csquotes}
\usepackage{algorithm} 
\usepackage{algpseudocode}
\usepackage{tikz}
\usetikzlibrary{arrows,shapes.geometric,positioning,calc}
\usepackage[colorinlistoftodos]{todonotes}
\usepackage{changepage}
\usepackage{subfig}
\usepackage{soul}
\usepackage{hyperref}
\usepackage{todonotes}
\usepackage{enumerate}

\newtheorem{thm}{Theorem}
\newtheorem{prop}[thm]{Proposition}
\newtheorem{cor}[thm]{Corollary}
\newtheorem{lem}[thm]{Lemma}

\newsavebox{\qedB}

\sbox{\qedB}{\setlength{\unitlength}{1mm}
 \begin{picture}(4,4)(0,0)
  \thinlines
  {\put(0,0){\framebox(2.83,2.83){}}}%
  {\put(1.17,1.17){\framebox(2.83,2.83){}}}%
  {\put(0,0){\framebox(4,4){}}}%
  {\put(1.17,1.17){{\rule{1ex}{1ex} }}}%
 \end{picture}}
 
\newcommand{\QEDB}{\ifmmode\def\next{\tag"\usebox{\qedB}"}%
 \else\let\next=\relax
 {\unskip\nobreak\hfil\penalty50
 \hskip2em\hbox{}\nobreak\hfil\usebox{\qedB}
 \parfillskip=0pt \finalhyphendemerits=0\penalty-100\bigskip}\fi\next}

\newcommand{\Alphabet}{\hbox{\rm Alph}}

\newcommand{\fac}{\hbox{\rm Fac}}

\newcommand{\Class}{\hbox{\rm Class}}
\newcommand{\Index}{\hbox{\rm Index}}

\newcommand{\prim}{\hbox{\rm Prim}}

\newcommand{\N}{\mathbb N}

\newcommand{\bprop}{\begin{prop}}
\newcommand{\eprop}{\end{prop}}
\newcommand{\bcor}{\begin{cor}}
\newcommand{\ecor}{\end{cor}}
\newcommand{\blem}{\begin{lem}}
\newcommand{\elem}{\end{lem}}

\title{On the number of squares in a finite word \thanks{with the support of NSERC (Canada)}}
\toctitle{Square complexity}

\authorrunning{S. Brlek, S. Li}
\tocauthor{Sre\v{c}ko Brlek, Shuo Li}
\author{S. Brlek\inst{1} \and S. Li\inst{1}}

\institute{Laboratoire de Combinatoire et d'Informatique Mathématique,\\
Université du  Québec \`a Montréal,\\
CP 8888 Succ. Centre-ville, Montréal (QC) Canada H3C 3P8\\
\email{brlek.srecko@uqam.ca, shuo.li@lacim.ca} 
}

\begin{document}

\maketitle

%%%%%%%%%%%%%%%%%%%%%%%%%%%%%%%%%%%%%%%%%%%%%%%%%%%%
% Sections
%%%%%%%%%%%%%%%%%%%%%%%%%%%%%%%%%%%%%%%%%%%%%%%%%%%%

\begin{abstract}
A {\em square} is a word of the form $uu$.
In this paper we prove that for a given finite word $w$, the number of distinct square factors of $w$ is bounded by $|w|-|\Alphabet(w)|+1$, where $|w|$ denotes the length of $w$ and $|\Alphabet(w)|$ denotes the number of distinct letters in $w$. This result answers a conjecture of Fraenkel and Simpson stated in 1998.
\end{abstract}

\section{Introduction}
The study of patterns in a word is 
related to several topics and it has
%a well studied topic
% for the numerous and
various applications in several research domains, ranging from practical applications to theoretical considerations.
Among these patterns, we can mention palindromes and repetitions such as squares, cubes, periods, overlaps, etc.
They appear in  different contexts in computer science such as data compression (LZW overlap), searching algorithms, structure of indexes \cite{lothaire3}, digital geometry \cite{bgl}, in number theory about Diophantine approximation and transcendence statements~\cite{BorisA}, and in physics in connection with Schr\"odinger operators~\cite{ABCD}, among several others.

Palindromic complexity has received noticeable attention and a comprehensive survey was  provided by Allouche et al.~\cite{ABCD}.
%The role of palindromic factors in words is ubiquitous.
%It gives an insight on the intrinsic structure, due to its connection with the usual complexity, the characterization it provides for Sturmian words~\cite{deLu} or the relation with the notion of recurrence~\cite {BrlekLad}.
%Many other connections exist and we refer to the survey~\cite{ABCD} for further reading. 
In particular, Droubay, Justin and Pirillo established that the number of distinct palindromic factors of a word is bounded by its length plus one and  that finite Sturmian (and even episturmian) words realize it~\cite{DrJuPi}.

%In this framework, Brlek and Reutenauer proposed the following identity~\cite{BrlekR}
%\begin{equation}
%\label{eq:pal}
 %2\D_p(w) = \sum_{n=0}^{|w|} C_w(n+1) - C_w(n) + 2 - P_w(n+1) - P_w(n)
%\end{equation}
%linking the palindromic defect $\D_p$ with the factor and palindromic complexities $C_w$ and $P_w$, motivated by the fact that it is (almost) obviously true for finite words. They provided several examples of infinite words satisfying it, including periodic ones, the Thue-Morse word, all Sturmian ones, and also the Oldenburger exponent trajectory. Conjectured  for languages closed by reversal, it  was established by Balkov\`a, Pelantov\`a and Starosta~\cite{BalkovaPS}.

%Later, the notion of defect was extended to $\sigma$-palindromes where $\sigma$ is an involution, also known as anti-palindromes~\cite{BrLaN}.
%Extending Identity~\eqref{eq:pal}, Brlek and Lafreni\`ere~\cite{nadia} proved the identity
%\begin{equation}
%\label{eq:spal}
 %2\D_\sigma(w) = \sum_{n=0}^{|w|} C_w(n+1)- C_w(n) + 2 - P_{\sigma,w}(n+1) - P_{\sigma,w}(n)
%\end{equation}
%which links the $\sigma$-palindromic defect $\D_\sigma$ with the factor and $\sigma$-palindromic complexities $C_w$ and $P_{\sigma,w}$,
%for finite words, periodic infinite words and infinite words whose language is closed under $\sigma$-reversal.
%This identity was presented in~\cite{SrBr2013Cetraro}.

%In this paper we consider first the language $L_S$ of words whose square factors belong to a fixed and finite set $S$ of squares.
For the square factors the situation seems to be significantly harder to handle.
Indeed, Fraenkel and Simpson conjectured in~\cite{FraenkelS98} that the number of distinct non-empty squares of a finite word $w$ is bounded by its length $|w|$ and they proved that this number is bounded by $2|w|$. 
%The best known upper bound is $\floor{11n/6}$ and was established by Deza \emph{et al.} in~\cite{DezaFT15}.
%, while for binary alphabets it was refined to ... \cite{JonoskaMS}. 
%Fraenkel and Simpson proved in~\cite{FraenkelS98} that
After that Ilie~\cite{Ilie} strengthened this bound to $2|w|-\Theta(\log(|w|))$; Lam~\cite{lam}
improved this result to $\frac{95}{48}|w|$; Deza, Franek and Thierry~\cite{DezaFT15}
achieved a bound of $\frac{11}{6}|w|$; Thierry~\cite{thie} refined this bound to $\frac{3}{2}|w|$. We remark that the Fraenkel and Simpson's conjecture implies that the square complexity of a finite word is also bounded by its length plus one by counting the empty word. In this article, we prove the following theorem:

\begin{thm}
\label{th:sw}
For any finite word $w$, let $S(w)$ denote the number of its distinct square factors, $|w|$ denote the length of $w$ and $|\Alphabet(w)|$ denote the number of distinct letters in $w$, then we have $$S(w) \leq |w|-|\Alphabet(w)|+1.$$
\end{thm}

This result is a stronger version of the conjecture of Fraenkel and Simpson stated in~\cite{FraenkelS98}.

The strategy of the proof is as follows: for a given word $w$, we first recall the definition of Rauzy graph and define small circuits in it; we then prove that the total number of small circuits in the union of the Rauzy graphs of $w$ is bounded by $|w|-|\Alphabet(w)|$; finally we conclude by proving that there exists an injection from the set of distinct nonempty squares of $w$ to the set of small circuits in the union of the Rauzy graphs of $w$.\\
%In this paper, we  propose the following identity
%\begin{equation}
%\label{eq:squares}
 %2\D_s(w) = \sum_{n=0}^{|w|} C_w(n+1) - C_w(n) + 2 - S_w(n+1) - S_w(n)
%\end{equation}
%that links the square defect $\D_s$ of a word with its factor and square complexities $C_w$ and $S_w$.
%We prove that this identity is verified for finite words as well as for some classes of infinite words.

%The paper is organized as follows.
%Section~\ref{sec:prel} contains preliminary definitions about words and squares.
%In Section~\ref{sec:few}, we consider the case of words having only a finite number of squares.
%In particular, we show that infinite languages having only a finite number of distinct squares are not context-free.
%The square defect is introduced in Section~\ref{sec:defect}, where it is analyzed in the case of finite words and infinite periodic words.
% and infinite aperiodic words obtained as fixed point of a primitive morphism.\todo{effacer}
%In Section~\ref{sec:br}, we introduce the generalized Brlek-Reutenauer identity for squares.
%We prove that this identity is satisfied by finite words as well as infinite periodic words.
% and a special case of infinite aperiodic words: strict standard episturmian words.
%Finally, in Section~\ref{sec:concl} we conclude with some remarks and open problems.

\section{Preliminaries}

Let us recall the basic terminology about words. By {\em word} we mean a finite concatenation of symbols $w = w_1 w_2 \cdots w_{n}$, with $n \in \N$. The {\em length} of $w$, denoted $|w|$, is $n$ and we say that the symbol $w_i$ is at the {\em position} $i$. The set $\Alphabet(w)=\left\{w_i| 1\leq i \leq n\right\}$ is called the {\em alphabet} of $w$ and its elements are called {\em letters}. Let $|\Alphabet(w)|$ denote the cardinality of $\Alphabet(w)$. A word of length $0$ is called the {\em empty word} and it is denoted by $\varepsilon$. For any word $u$, we have $u=\varepsilon u=u\varepsilon$.

A word $u$ is called a {\em factor} of $w$ if $w = pus$ for some words $p,s$. The set of all factors of $w$ is denoted by $\fac(w)$. For any integer $i$ satisfying $1 \leq i \leq |w|$, let $L_w(i)$ be the set of all length-$i$ factors of $w$ and let $C_w(i)$ denote the cardinality of $L_i(w)$.

Two finite words $u$ and $v$ are {\em conjugate} when there exist words $x,y$ such that $u=xy$ and $v=yx$.
The conjugacy class of a word $w$ is denoted by $[w]$. Let $w=w_1w_2...w_t$ be a finite word, for any integer $i$ satisfying $1 < i \leq t$, let us define $w_{s}(i)=w_iw_{i+1}...w_{t}$ and $w_p(i)=w_1w_2...w_{i-1}$ and let $w_s(1)=w$ and $w_p(1)=\epsilon$. 
%For all $i$ such that $i > |v|$, let us define $v_s(i)=v_s(j)$ (resp. $v_p(i)=v_p(j)$) if $i \equiv j \; (\mod |v|)$. 
Thus, $[w]=\left\{w(i)|w(i)=w_s(i)w_p(i), i=1,2,...t\right\}$.

For any natural number $k$, we define the {\em $k$-power} of a finite word $u$ to be the concatenation of $k$ copies of $u$, and it is denoted by $u^k$. Particularly, a \emph{square} is a word $w$ of the form $w=uu$. 
A word $w$ is said to be {\em primitive} if it is not a power of another word. Let $\prim(w)$ denote the set of primitive factors of $w$. For any word $u$ and any rational number $\alpha$, the $\alpha$-power of $u$ is defined to be $u^au_0$ where $u_0$ is a prefix of $u$, $a$ is the integer part of $\alpha$, and $|u^au_0|=\alpha |u|$. The $\alpha$-power of $u$ is denoted by $u^{\alpha}$. Let $\alpha$ be a positive rational number larger than $1$, the word $w$ is said to be {\em of the period $\alpha$} if there exists a word $u$ such that $w=u^{\alpha}$.\\

Here we recall some basic lemmas concerning repetitions:\\

\begin{lemma}[Fine and Wilf~\cite{fiwi}]
\label{period}
Let $w$ be a word having $k$ and $l$ for periods. If $|w| \geq k+l-\gcd (k,l)$ then $\gcd(k,l)$ is also a period of $w$.
\end{lemma}

\begin{lemma}[Lyndon and Schützenberger~\cite{lyndon}]
\label{two-words}
Let $x$ and $y$ be two words such that $xy = yx$. Then there exist a primitive $p$ and two non-negative integers $i, j$ such that $x = p^i$ and y = $p^j$.
\end{lemma}

\begin{lemma}[Lyndon and Schützenberger~\cite{lyndon}]
\label{three-words}
Let $x,y,z$ be three words such that $xy = yz$ and $|y| \neq 0$. Then there exist two words $u, v$ with $|u|\neq 0$ and some positive integer $i$ such that 
$x = uv, y = (uv)^iu, z = vu$.
\end{lemma}

Here we recall some elementary definitions and proprieties concerning graphs from Berge~\cite{berge}.\\

A {\em (non-oriented) graph} consists of a nonempty set of {\em vertices} $V$ and a set of {\em edges} $E$. A vertex $a$ represents an endpoint of an edge and an edge joins two vertices $a, b$.
A {\em chain} is a sequence of edges $e_1,e_2, \cdots, e_k$, in which each edge $e_i$ ($1<i<k$) has one vertex in common with $e_{i-1}$ and the other vertex in common with $e_{i+1}$. A {\em cycle} is a finite chain which begins with a vertex $x$ and ends at the same vertex. A graph is called {\em connected} if for any couple of vertices $a,b$ in this graph, there exists a chain which begins at $a$ and ends at $b$.

A graph is called {\em oriented} if its edges are oriented from one endpoint to the other. An oriented graph is called weakly connected if it is connected as a non-oriented graph.

Let $G$ be a connected graph and let $\left\{e_1,e_2 \cdots e_l\right\}$, $\left\{v_1,v_2 \cdots v_s\right\}$ denote respectively the edge set and the vertex set of $G$. The number $\chi(G)=l-s+1$ is called the {\em cyclomatic number} of $G$.

Let $C$ be a cycle in $G$. A vector $\mu(C)=(c_1,c_2 \cdots c_l)$ in the $l$-dimensional space $\mathbb{R}^l$ is called the {\em vector-cycle corresponding to $C$} if $c_i$ is the number of visits of the edge $e_i$ in the cycle $C$ for all $i$ satisfying $1 \leq i \leq l$. The cycle $C_1,C_2, \cdots, C_k,...$ are said to be {\em independent} if their corresponding vectors are linearly independent. 

\begin{lemma}[Theorem 2, Chapter 4 in~\cite{berge}]
\label{book}
the cyclomatic number of a graph is the maximum number of independent cycles in this graph.
\end{lemma}

\section{Rauzy graph and small circuits}
We first recall the definition of Rauzy graph. Let $w$ be a word of length $k$. For any integer $n$ such that $1\leq n\leq k$, let the Rauzy graph $\Gamma_n(w)$ be an oriented graph whose set of vertices is $L_w(n)$ and the set of edges is $L_{w}(n+1)$;
an edge $e \in L_{w}(n+1)$ starts at the vertex $u$ and ends at the vertex $v$, if $u$ is a
prefix and $v$ is a suffix of $e$. We remark that $\Gamma_n(w)$ is a weakly connected graph for any $n$. Let us define $\Gamma(w)=\cup_{n=1}^{k}\Gamma_n(w)$.

Let $\Gamma_i(w)$ be a Rauzy graph of $w$. A sub-graph in $\Gamma_i(w)$ is called an {\em elementary circuit} if there are $j$ distinct vertices $v_1,v_2,..., v_j$ and $j$ distinct edges $e_1,e_2,...,e_j$ for some integer $j$, such that for each $k$ with $1 \leq k \leq j-1$, the edge $e_k$ starts at $v_k$ and ends at $v_{k+1}$, and for the edge $e_j$, it starts at $v_j$ and ends at $v_1$, further, $j$ is called the {\em size} of the circuit. We let this circuit be denoted by $\left\{\left\{v_1,v_2,..., v_j\right\},\left\{e_1,e_2,...,e_j\right\}\right\}$. The {\em small circuits} in the graph $\Gamma_i(w)$ are those elementary circuits whose sizes are no larger than $i$. 

\begin{lemma}
\label{small-circuit}
Let $w$ be a finite word and let $\Gamma_r(w)$ be a Rauzy graph of $w$ for some $r$. Then for any small circuit $C$ on $\Gamma_r(w)$, there exists a unique primitive word $q$, up to conjugacy, such that $|q| \leq r$ and the vertex set of $C$ is $\left\{p^{\frac{r}{|p|}}| p \in [q]\right\}$ and its edge set is $\left\{p^{\frac{r+1}{|p|}}| p \in [q]\right\}$.
\end{lemma}

\begin{proof}
Let $C$ be a small circuit on $\Gamma_r(w)$ of size $l$, then it contains $l$ distinct vertices, say $v_1,v_2,...,v_l$, and $l$ distinct edges, say $e_1,e_2,...,e_l$. For each $i$ satisfying $1 \leq i \leq l$, let us define $p_i$ to be a word by concatenating consecutively the last letter of words $e_i,e_{i+1},...,e_{i+l-1}$ with $e_{r}=e_{r-l}$ if $r\geq l+1$. The words $p_1,p_2,...,p_l$ are pairwisely conjugate. Further, for each $i$ satisfying $1 \leq i \leq l$, from the fact that the edges in the order of $e_i,e_{i+1},...,e_{i+l-1}$ form a circuit, we can deduce that $v_i$ is a suffix of $v_ip_i$. Thus, there exists a word $p'_i$ of length $l$ such that $v_ip_i=p'_iv_i$. From Lemma~\ref{three-words}, there exists two words $u_i,t_i$ with $|u_i| \neq 0$ and a nonnegative integer $j$ such that $p'_i=u_it_i, p_i=t_iu_i$ and $v_i=(u_it_i)^ju_i$. Consequently, $v_i=(p'_i)^{\frac{r}{l}}$. Further, as $p'_i$ and $p_i$ are conjugate, the words $p'_1,p'_2,...,p'_l$ are also pairwisely conjugate. Thus, the vertex set of $C$ is $\left\{p^{\frac{r}{|p|}}| p \in [p_1]\right\}$ and its edge set is $\left\{p^{\frac{r+1}{|p|}}| p \in [p_1]\right\}$.

Here we prove the primitivity of $p_1$. If $p_1$ is not primitive, then there are two distinct integers $r,s$ satisfying $1 \leq r,s \leq l$ such that $p'_r=p'_s$, and further, $v_r=v_s$. This contradicts the hypothesis that $C$ contains $l$ distinct vertices.

\end{proof}

From the previous lemma, each small circuit can be uniquely identified by two parameters: a primitive word $p$ defined in the lemma and an integer $i$ such that $\Gamma_w(i)$ is the Rauzy graph in which the circuit is located. Let all small circuits be denoted by $C(p,i)$ with the parameters defined as above.

\begin{example}
Let us consider the word $u=aababa$, here we draw $\Gamma_u(1), \Gamma_u(2),\Gamma_u(3)$ as follows:
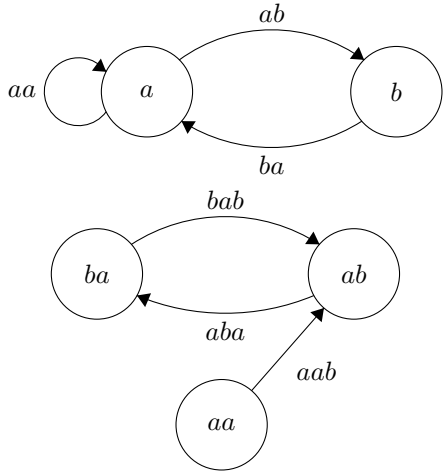
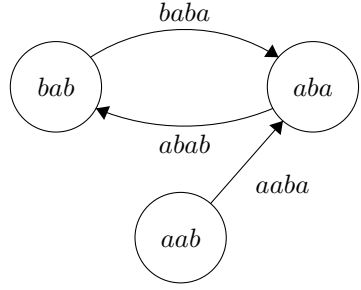
\begin{center}
\begin{tikzpicture}[scale=0.2]
\tikzstyle{every node}+=[inner sep=0pt]
\draw [black] (25.6,-6.5) circle (3);
\draw (25.6,-6.5) node {$b$};
\draw [black] (9.2,-6.5) circle (3);
\draw (9.2,-6.5) node {$a$};
\draw [black] (11.322,-4.394) arc (126.39894:53.60106:10.243);
\fill [black] (23.48,-4.39) -- (23.13,-3.52) -- (22.54,-4.32);
\draw (17.4,-1.9) node [above] {$ab$};
\draw [black] (23.333,-8.451) arc (-57.14522:-122.85478:10.936);
\fill [black] (11.47,-8.45) -- (11.87,-9.3) -- (12.41,-8.46);
\draw (17.4,-10.7) node [below] {$ba$};
\draw [black] (6.52,-7.823) arc (324:36:2.25);
\draw (1.95,-6.5) node [left] {$aa$};
\fill [black] (6.52,-5.18) -- (6.17,-4.3) -- (5.58,-5.11);
\end{tikzpicture}
\end{center}

\begin{center}
\begin{tikzpicture}[scale=0.2]
\tikzstyle{every node}+=[inner sep=0pt]
\draw [black] (20.1,-6.2) circle (3);
\draw (20.1,-6.2) node {$ab$};
\draw [black] (3.2,-6.2) circle (3);
\draw (3.2,-6.2) node {$ba$};
\draw [black] (11.4,-16.2) circle (3);
\draw (11.4,-16.2) node {$aa$};
\draw [black] (13.37,-13.94) -- (18.13,-8.46);
\fill [black] (18.13,-8.46) -- (17.23,-8.74) -- (17.98,-9.4);
\draw (16.29,-12.65) node [right] {$aab$};
\draw [black] (5.463,-4.244) arc (123.23325:56.76675:11.29);
\fill [black] (17.84,-4.24) -- (17.44,-3.39) -- (16.89,-4.22);
\draw (11.65,-1.9) node [above] {$bab$};
\draw [black] (17.469,-7.632) arc (-67.17079:-112.82921:14.999);
\fill [black] (5.83,-7.63) -- (6.37,-8.4) -- (6.76,-7.48);
\draw (11.65,-9.31) node [below] {$aba$};
\end{tikzpicture}
\end{center}

\begin{center}
\begin{tikzpicture}[scale=0.2]
\tikzstyle{every node}+=[inner sep=0pt]
\draw [black] (20.1,-6.2) circle (3);
\draw (20.1,-6.2) node {$aba$};
\draw [black] (3.2,-6.2) circle (3);
\draw (3.2,-6.2) node {$bab$};
\draw [black] (11.4,-16.2) circle (3);
\draw (11.4,-16.2) node {$aab$};
\draw [black] (13.37,-13.94) -- (18.13,-8.46);
\fill [black] (18.13,-8.46) -- (17.23,-8.74) -- (17.98,-9.4);
\draw (16.29,-12.65) node [right] {$aaba$};
\draw [black] (5.463,-4.244) arc (123.23325:56.76675:11.29);
\fill [black] (17.84,-4.24) -- (17.44,-3.39) -- (16.89,-4.22);
\draw (11.65,-1.9) node [above] {$baba$};
\draw [black] (17.469,-7.632) arc (-67.17079:-112.82921:14.999);
\fill [black] (5.83,-7.63) -- (6.37,-8.4) -- (6.76,-7.48);
\draw (11.65,-9.31) node [below] {$abab$};
\end{tikzpicture}
\end{center}
There are three small circuits in these graphs: $\left\{\left\{a\right\},\left\{aa\right\}\right\}$, $\left\{\left\{ab,ba\right\},\left\{aba,bab\right\}\right\}$ and $\left\{\left\{aba,bab\right\},\left\{abab,baba\right\}\right\}$. With the notation defined as above, these circuits are denoted respectively by $C(a,1), C(ab,2)$ and $C(ab,3)$. We remark that, in $\Gamma_1(u)$, there exists a circuit $\left\{\left\{a,b\right\},\left\{ab,ba\right\}\right\}$ which is not small.

\end{example}

Let $\prec$ be a lexicographic order over $\fac(w)$, for a given small circuit, let us define the {\em maximal edge} of this circuit to be the lexicographically greatest element in its edge set.

\begin{lemma}
\label{common}
Each small circuit contains exactly one maximal edge. Further, the maximal edges of the small circuits in the same Rauzy graph are pairwisely distinct. 
\end{lemma}

\begin{proof}

From Lemma~\ref{small-circuit} the edge set of any small circuit is of the form $\left\{p^{t}| p \in [q]\right\}$ such that $t > 1$ and $q$ is primitive. Thus, $u^t$ is the lexicographically greatest element in $\left\{p^{t}| p \in [q]\right\}$ if and only if $u$ is the  lexicographically greatest element in $[q]$. Hence the maximal edge of each small circuit is unique.\\

For the second part, let us suppose that $v$ is the maximal edge of two distinct small circuits $C(p,r)$ and $C(q,r)$. From Lemma~\ref{small-circuit}, there are two positive integers $i\geq 1,j \geq 1$ and two words $p',q'$ which are respectively a prefix of $p, q$ such that $v=p^ip'=q^jq'$,where $p',q'$ can be $\varepsilon$. First we state that $|p'| \neq |q'|$. Otherwise, $p^i=q^j$. However, in this case, from Lemma~\ref{period}, $p$ and $q$ can not be both primitive. Without loss of generality, let us suppose $|p'|>|q'|$, thus $|p^i|<|q^j|$. Then there exists a prefix $p''$ of $p$ such that $|p''| <|p|$ and $q^j=p^ip''$. However, $p''$ is also a prefix of $p^i$. Thus, there exists a word $w$ such that $q^j=p''wp''$ and $p^i=p''w$.

Now let us prove that $v$ cannot be the maximal edge of $C(q,r)$. It is enough to prove that $q^j=p''wp''$ is not the lexicographically greatest element in $[q^j]$. Let us consider two other conjugates of $q^j$: $(p'')^2w$ and $w(p'')^2$. We claim that $(p'')^2w \neq p''wp''$. Otherwise, from Lemma~\ref{two-words}, there exists a primitive $a$ and two positive integers $t_2>t_1>0$ such that $p''=a^{t_1}, p''w=a^{t_2}$. However, in this case, $p^j=p''w=a^{t_2}$ and $|a| < |p|$. Once more from Lemma~\ref{period}, $p$ and $a$ can not be both primitive. If $(p'')^2w \prec p''wp''$, then $p''w \prec wp''$, but if so, $p''wp'' \prec w(p'')^2$. We conclude. 
\end{proof}

Now let us define the {\em Circuit Arrangement Order} $\prec_{CAO}$ for the small circuits on $\Gamma(w)$. Let $C_1$ and $C_2$ be two small circuits in a same Rauzy graph $\Gamma_r(w)$, we define $C_1 \prec_{CAO} C_2$ if $e_1 \prec e_2$, where $e_i$ is the maximal edge of the cycle $C_i$. From Lemma \ref{common}, the small circuits in the same Rauzy graph are well-ordered under $\prec_{CAO}$. 

\begin{lemma}
\label{independent}
Let $w$ be a finite word and let $\Gamma_r(w)$ be a Rauzy graph of $w$ for some $r$. Then the small circuits in $\Gamma_r(w)$ are independent.
\end{lemma}

\begin{proof}
Let $C_1,C_2, \cdots, C_k$ be the small circuits in $\Gamma_r(w)$ satisfying $$C_1 \prec_{CAO} C_2 \prec_{CAO} \cdots \prec_{CAO}C_k.$$ If there exist some real numbers $\alpha_1,\alpha_2, \cdots, \alpha_k$ such that $\sum_{i=1}^k \alpha_i\mu(C_i)=0$, where $\mu(C_i)$ is the vector-cycle corresponding to $C_i$, then we claim that $\alpha_k=0$. In fact, if we let $e_{(C_k)}$ be the maximal edge of $C_k$, then from Lemma~\ref{common}, the coefficient $c_{(C_k)}$ in $\mu(C_i)$ is $1$ if $i=k$ and $0$ if $i < k$. Thus, $\alpha_k=0$. Doing this argument recurrently, we prove that $\alpha_i=0$ for all $i$. Consequently, $C_1,C_2, \cdots, C_k$ are independent.
\end{proof}

\begin{example}
Let us consider the word $v=abaaabaabaaaba$, the $\Gamma_5(v)$ is drown as follows:
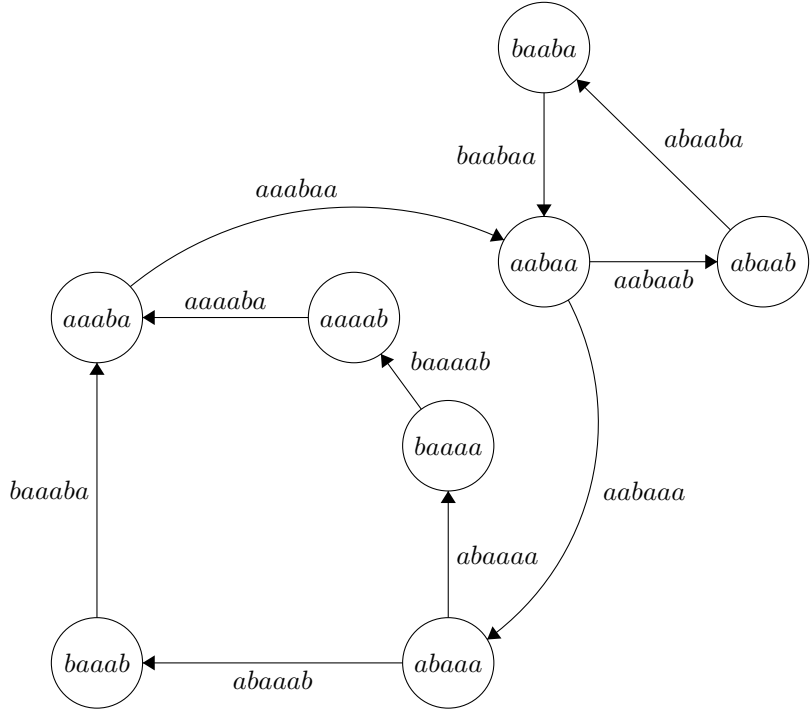
\begin{center}
\begin{tikzpicture}[scale=0.2]
\tikzstyle{every node}+=[inner sep=0pt]
\draw [black] (6.3,-21.1) circle (3);
\draw (6.3,-21.1) node {$aaaba$};
\draw [black] (35.7,-17.4) circle (3);
\draw (35.7,-17.4) node {$aabaa$};
\draw [black] (6.3,-44) circle (3);
\draw (6.3,-44) node {$baaab$};
\draw [black] (29.4,-44) circle (3);
\draw (29.4,-44) node {$abaaa$};
\draw [black] (23.2,-21.1) circle (3);
\draw (23.2,-21.1) node {$aaaab$};
\draw [black] (29.4,-29.6) circle (3);
\draw (29.4,-29.6) node {$baaaa$};
\draw [black] (35.7,-3.2) circle (3);
\draw (35.7,-3.2) node {$baaba$};
\draw [black] (50.1,-17.4) circle (3);
\draw (50.1,-17.4) node {$abaab$};
\draw [black] (20.2,-21.1) -- (9.3,-21.1);
\fill [black] (9.3,-21.1) -- (10.1,-21.6) -- (10.1,-20.6);
\draw (14.75,-20.6) node [above] {$aaaaba$};
\draw [black] (27.63,-27.18) -- (24.97,-23.52);
\fill [black] (24.97,-23.52) -- (25.04,-24.46) -- (25.84,-23.88);
\draw (26.88,-23.97) node [right] {$baaaab$};
\draw [black] (29.4,-41) -- (29.4,-32.6);
\fill [black] (29.4,-32.6) -- (28.9,-33.4) -- (29.9,-33.4);
\draw (29.9,-36.8) node [right] {$abaaaa$};
\draw [black] (8.49,-19.053) arc (129.36441:64.98155:23.251);
\fill [black] (33.07,-15.96) -- (32.56,-15.17) -- (32.13,-16.07);
\draw (19.43,-13.17) node [above] {$aaabaa$};
\draw [black] (37.293,-19.938) arc (27.28469:-53.93376:17.769);
\fill [black] (31.96,-42.45) -- (32.9,-42.38) -- (32.31,-41.57);
\draw (39.55,-32.6) node [right] {$aabaaa$};
\draw [black] (26.4,-44) -- (9.3,-44);
\fill [black] (9.3,-44) -- (10.1,-44.5) -- (10.1,-43.5);
\draw (17.85,-44.5) node [below] {$abaaab$};
\draw [black] (6.3,-41) -- (6.3,-24.1);
\fill [black] (6.3,-24.1) -- (5.8,-24.9) -- (6.8,-24.9);
\draw (5.8,-32.55) node [left] {$baaaba$};
\draw [black] (38.7,-17.4) -- (47.1,-17.4);
\fill [black] (47.1,-17.4) -- (46.3,-16.9) -- (46.3,-17.9);
\draw (42.9,-17.9) node [below] {$aabaab$};
\draw [black] (47.96,-15.29) -- (37.84,-5.31);
\fill [black] (37.84,-5.31) -- (38.05,-6.22) -- (38.76,-5.51);
\draw (46.2,-9.82) node [above] {$abaaba$};
\draw [black] (35.7,-6.2) -- (35.7,-14.4);
\fill [black] (35.7,-14.4) -- (36.2,-13.6) -- (35.2,-13.6);
\draw (35.2,-10.3) node [left] {$baabaa$};
\end{tikzpicture}
\end{center}

In this graph, there exist three small circuits: $C(aaaab,5), C(aaab,5), C(aab,5)$. If we define $b \prec a$, the maximal edges of these circuits are respectively $aaaaba$, $aaabaa$ and $aabaab$. We can check that $aaaaba$ is only in $C(aaaab,5)$ and $aabaab$ is only in $C(aab,5)$, thus, the three circuits in the graph are independent.
\end{example}

\begin{lemma}
\label{bound}
Let $w$ be a finite word, then for any Rauzy graph $\Gamma_r(w)$, the number of small circuits in this graph, which is denoted by $sc_r$, satisfies that $$sc_r \leq C_w(r+1)-C_w(r)+1.$$
\end{lemma}

\begin{proof}
It is a direct consequence of Lemma~\ref{book}, Lemma~\ref{independent} together with the facts that $\Gamma_r(w)$ is weakly-connected and the number of edges and vertices in $\Gamma_r(w)$ are respectively $C_w(r+1)$ and $C_w(r)$. 
\end{proof}

\begin{lemma}
\label{total-number}
Let $w$ be a finite word, then the total number of small circuits on $\Gamma(w)$ is bounded by $|w|-|\Alphabet(w)|$.
\end{lemma}

\begin{proof}
Let $sc(w)$ denote the total number of small circuits on $\Gamma(w)$, then from the previous lemma, we have:
\begin{align*}
sc(w)&= \sum_{i=1}^{|w|} sc_i\\
&\leq \sum_{i=1}^{|w|}C_w(i+1)-C_w(i)+1\\
&\leq |w|+C_{w}(|w|+1)-C_{w}(1)\\
&\leq |w|-|\Alphabet(w)|.\\
\end{align*}

\end{proof}

\section{Injection from nonempty squares to small circuits}

Let $w$ be a finite word. For any square $u^2 \in \fac(w)$, there exists a primitive word $v \in \prim(w)$ and a positive integer $n$ such that $u^2=v^{2n}$. Thus, we can gather the squares in terms of their smallest period. Let  $v \in \prim(w)$, a square $u^2 \in \fac(w)$ is called {\em in the class of $v$} if there is a primitive word $t \in [v]$ and a positive integer $n$ such that $u^2=t^{2n}$. Let $\Class(v)$ be the set of all squares in the class of $v$. Two classes $\Class(u)$ and $\Class(v)$ are equal only if $u$ and $v$ are conjugate.

Now given a class $\Class(v)$, let us define its {\em index} to be an integer $\Index(v)$ such that  $\Index(v)=\max\left\{n| n \in \mathbb{N^+}, \exists u \in [v],\; \textrm{such that}\; u^{2n} \in \fac(w)\right\}$. Without loss of generality, a class is denoted by $\Class(v)$ if $v$ is a primitive such that $v^{2n}$ is in the class where $n$ is the index of this class. In other words, $v^{2\Index(v)} \in \Class(v)$. From the definition, the elements in $\Class(v)$ are all of the form $v_s(i)v^{2j-1}v_p(i)$ with $1\leq i \leq |v|$ and $1\leq j \leq \Index(v)$, let $v(i,j)$ denote the square $v_s(i)v^{2j-1}v_p(i)$. 

\begin{example}
Let $r=abcabcabcabca$. The squares $(abc)^2, (bca)^2,(cab)^2,(abc)^4, (bca)^4$ are in the same class. The index of this class is $4$ and this class can be denoted by $\Class(abc)$ or $\Class(bca)$. However, it cannot be denoted by $\Class(cab)$, because $(cab)^4 \not \in \fac(r)$.
\end{example}

\begin{lemma}
Let $v$ be a primitive factor of a finite word $w$ such that $|v|=l$ and that $\Class(v)$ is nonempty. If $\Index(v) \geq 2$, then for any couple of integers $(i,j)$ satisfying $1\leq i \leq l$ and $1\leq j \leq \Index(v)$, there exists a small circuit $C(v, jl+i-1)$ in $\Gamma_{jl+i-1}(w)$. Further, there exists an injective function $f_v$ which associates each square $v(i,j)$ in $\Class(v)$ to the small circuit $C(v, jl+i-1)$ in $\left\{C(v, jl+i-1)| 1 \leq i \leq l, 1\leq j \leq \Index(v)\right\}$.
\end{lemma}

\begin{proof}
For a given couple of integers $(i,j)$, it is enough to prove that 
$$\left\{u^{j+\frac{i}{l}}|u \in [v]\right\} \subset \fac(w).$$ If it is the case, then there exists a circuit in $\Gamma_w(jl+i-1)$ such that its edge set is $\left\{u^{j+\frac{i}{l}}|u \in [v]\right\}$, further, it can be identified by $C(v, jl+i-1)$. From the hypothesis that $\Index(v) \geq 2$ and $v^{2\Index(v)} \in \fac(w)$, we have $v^{j+2} \in \fac(w)$ for all $j$ satisfying $1\leq j \leq \Index(v)$ and easily check that $\left\{u^{j+\frac{i}{l}}|u \in [v]\right\} \subset \fac(v^{j+2})$.

The injectivity of the function is from the fact that the circuits $C(v, jl+i-1)$ are pairwisely distinct. In fact, if there are two pairs of integers $(j_1,i_1), (j_2,i_2)$ satisfying $1 \leq i_1,i_2 \leq l$, $j_1 \geq j_2$ and $j_1l+i_1-1=j_2l+i_2-1$, then $(j_1-j_2)l=i_1-i_2$. However, $i_1-i_2 <l$, hence $0\leq j_1-j_2<1$. Thus, $j_1=j_2$ and $i_1=i_2$.
\end{proof}

\begin{lemma}
Let $v$ be a primitive factor of a finite word $w$ such that $|v|=l$ and that the size of $\Class(v)$ is $t$. If $\Index(v)=1$, then for each integer $i$ such that $1 \leq i \leq t$, there exists a small circuit $C(v, l+i-1)$ in $\Gamma_{l+i-1}(w)$. Thus, there exists a bijective function $f_v$ which associates each square in $\Class(v)$ to one of the small circuits in $\left\{C(v, l+i-1)| 1 \leq i \leq t\right\}$.
\end{lemma}

\begin{proof}
It is enough to prove that for any positive integer $i$ satisfying $1 \leq i \leq t$, $$S_i=\left\{u^{1+\frac{i}{l}}|u \in [v]\right\} \subset \fac(w).$$ If so, there exists a circuit in $\Gamma_{l+i-1}$ such that its edge set is exactly $S_i$.  Let us first consider the elements in $\Class(v)$. Under the hypothesis that $\Index(v)=1$, the elements in $\Class(v)$ are all conjugates of $v^2$, thus they are all of the form $v_s(j)vv_p(j)$ for some $j$. Consequently, the size of $\Class(v)$ is no larger than the number of conjugates of $v$, thus, $t \leq l$. Now let us consider the elements in $S_i$. For any given $i$ satisfying $1 \leq i \leq t$ and any $u^{1+\frac{i}{l}} \in S_i$, there exists a couple of integers $(r,s)$ satisfying $s \equiv r+i \; (\mod l)$ and $u^{1+\frac{i}{l}}=v_s(r)v_p(s)$. For any $v_s(j)vv_p(j)\in \Class(v)$, a word $v_s(r)v_p(s) \in S_i$ is not a factor of $v_s(j)vv_p(j)$ if and only if $s <j <r$. Further, from the relation $s \equiv r+i \; (\mod l)$, we have $r-s=i$. Thus, for any $v_s(r)v_p(s) \in S_i$, there are at most $i-1$ distinct words $v_s(j)vv_p(j)\in \Class(v)$ which do not contain $v_s(r)v_p(s)$ as a factor. However, from the hypothesis, there are exactly $t \geq i$ elements in $\Class(v)$. Hence there exists at least one word in $\Class(v)$ such that $v_s(r)v_p(s)$ is a factor of this word. Thus, all elements in $S_i$ are factors of $w$.
\end{proof}

\begin{lemma}
\label{injective}
There exists an injective function $f$ from the set of nonempty squares of $w$ to the set of small circuits in $\Gamma(w)$.
\end{lemma}

\begin{proof}
Let us consider the function $f$ satisfying that for any $\Class(v)$, 
 \[f|_{\Class(v)}=f_v,\]
with $f_v$ defined in the previous lemmas. Here we prove that this function is injective. If $C(p,i)$ and $C(q,j)$ are two images of two squares under $f$ such that $C(p,i)=C(q,j)$, from Lemma~\ref{small-circuit}, $p=q$ and $i=j$. However, if $p=q$, the pre-image of these two circuits are in the same class, and from the previous lemmas, the function $f$ over a class is injective. Thus, these two squares are the same. Hence $f$ is injective.
\end{proof}

\begin{proof}[ of Theorem~\ref{th:sw}]
From Lemma~\ref{injective} and Lemma~\ref{total-number}, the number of distinct nonempty squares, $S(w)-1$, satisfies the equation
$$S(w)-1\leq sc(w) \leq |w|-|\Alphabet(w)|.$$
\end{proof}

\begin{example}
Let us consider the word $w=baababaababbbabbabbbab$. We can check that $|w|=22$ and there are 14 squares in $w$: 

$\varepsilon$, $ aa$, $bb$, $abab$, $baba$, $abaaba$, $bbabba$, $babbab$, $abbabb$, $babbbabb$, 

$bbabbbab$, $baababaaba$, $aababaabab$, $babbbabbabbbab$. 

\noindent
The nonempty squares are distributed in $8$ classes. In the following table we list these classes and their sizes:

$$ \begin{array}{|c|c|c|c|c|c|c|c|c|}
  \hline
  Class & \;\;\;a\;\;\; & \;\;\;b\;\;\; & \;\;ab\;\;\; & \;\;aba\;\; & \;\;abb\;\; & \;babb\;\; & \;baaba\; & babbbab \\
  \hline
  size & 1 & 1 & 2 & 1 & 3 & 2 & 2 & 1 \\
  \hline
 \end{array}$$
The images of these squares under the injection $f$ are listed as follows:\\

$ C(a,1), C(b,1), C(ab,2), C(ab,3), C(aba,3), C(abb,3), C(abb,4), $

$C(abb,5), C(babb,4), C(babb,5), C(baaba,5), C(baaba,6), C(babbbab,7). $

\noindent The Rauzy graphs of $w$ from $\Gamma_w(1)$ to $\Gamma_w(7)$ are provided in Appendix.

\end{example}

%%%%%%%%%%%%%%%%%%%%%%%%%%%%%%%%%%%%%%%%%%%%%%%%%%%%
% Biblio
%%%%%%%%%%%%%%%%%%%%%%%%%%%%%%%%%%%%%%%%%%%%%%%%%%%%

\bibliographystyle{splncs03}
\bibliography{biblio}
\newpage
\section{Appendix : graphs of $w=baababaababbbabbabbbab$}

Here we draw the Rauzy graphs of $w$ from $\Gamma_w(1)$ to $\Gamma_w(7)$, they are all the Rauzy graphs containing small cycles.

\begin{center}
\begin{tikzpicture}[scale=0.2]
\tikzstyle{every node}+=[inner sep=0pt]
\draw [black] (25.6,-6.5) circle (3);
\draw (25.6,-6.5) node {$b$};
\draw [black] (9.2,-6.5) circle (3);
\draw (9.2,-6.5) node {$a$};
\draw [black] (11.322,-4.394) arc (126.39894:53.60106:10.243);
\fill [black] (23.48,-4.39) -- (23.13,-3.52) -- (22.54,-4.32);
\draw (17.4,-1.9) node [above] {$ab$};
\draw [black] (23.333,-8.451) arc (-57.14522:-122.85478:10.936);
\fill [black] (11.47,-8.45) -- (11.87,-9.3) -- (12.41,-8.46);
\draw (17.4,-10.7) node [below] {$ba$};
\draw [black] (6.52,-7.823) arc (324:36:2.25);
\draw (1.95,-6.5) node [left] {$aa$};
\fill [black] (6.52,-5.18) -- (6.17,-4.3) -- (5.58,-5.11);
\end{tikzpicture}
\end{center}

\begin{center}
\begin{tikzpicture}[scale=0.2]
\tikzstyle{every node}+=[inner sep=0pt]
\draw [black] (21.7,-5.7) circle (3);
\draw (21.7,-5.7) node {$ab$};
\draw [black] (3.2,-5.7) circle (3);
\draw (3.2,-5.7) node {$aa$};
\draw [black] (3.2,-22.6) circle (3);
\draw (3.2,-22.6) node {$ba$};
\draw [black] (20.1,-23.9) circle (3);
\draw (20.1,-23.9) node {$bb$};
\draw [black] (5.688,-4.034) arc (117.86407:62.13593:14.467);
\fill [black] (19.21,-4.03) -- (18.74,-3.22) -- (18.27,-4.1);
\draw (12.45,-1.86) node [above] {$aab$};
\draw [black] (20.89,-8.586) arc (-19.70593:-75.46983:21.349);
\fill [black] (6.15,-22.05) -- (7.05,-22.34) -- (6.8,-21.37);
\draw (17.09,-17.64) node [below] {$aba$};
\draw [black] (3.2,-19.6) -- (3.2,-8.7);
\fill [black] (3.2,-8.7) -- (2.7,-9.5) -- (3.7,-9.5);
\draw (3.7,-14.15) node [right] {$baa$};
\draw [black] (4.136,-19.752) arc (158.0576:106.76664:22.913);
\fill [black] (18.78,-6.38) -- (17.87,-6.13) -- (18.16,-7.08);
\draw (7.98,-10.91) node [above] {$bab$};
\draw [black] (23.73,-7.896) arc (34.99474:-45.04289:11.08);
\fill [black] (22.48,-22.09) -- (23.4,-21.88) -- (22.7,-21.17);
\draw (26.32,-15.29) node [right] {$abb$};
\draw [black] (23.04,-23.365) arc (128.0546:-159.9454:2.25);
\draw (27.36,-26.32) node [right] {$bbb$};
\fill [black] (22.31,-25.91) -- (22.41,-26.85) -- (23.2,-26.23);
\end{tikzpicture}
\end{center}

\begin{center}
\begin{tikzpicture}[scale=0.2]
\tikzstyle{every node}+=[inner sep=0pt]
\draw [black] (20.9,-12.1) circle (3);
\draw (20.9,-12.1) node {$aba$};
\draw [black] (3.2,-3.2) circle (3);
\draw (3.2,-3.2) node {$aab$};
\draw [black] (3.2,-21.2) circle (3);
\draw (3.2,-21.2) node {$baa$};
\draw [black] (31.3,-12.1) circle (3);
\draw (31.3,-12.1) node {$bab$};
\draw [black] (42.8,-12.1) circle (3);
\draw (42.8,-12.1) node {$abb$};
\draw [black] (56.4,-4.2) circle (3);
\draw (56.4,-4.2) node {$bbb$};
\draw [black] (55.8,-21.2) circle (3);
\draw (55.8,-21.2) node {$bba$};
\draw [black] (6.182,-2.914) arc (90.20368:36.40755:16.291);
\fill [black] (19.35,-9.54) -- (19.28,-8.6) -- (18.47,-9.19);
\draw (15.87,-4.14) node [above] {$aaba$};
\draw [black] (19.661,-14.826) arc (-30.47953:-95.10288:14.224);
\fill [black] (6.14,-21.78) -- (6.89,-22.35) -- (6.98,-21.35);
\draw (16.22,-20.77) node [below] {$abaa$};
\draw [black] (3.2,-18.2) -- (3.2,-6.2);
\fill [black] (3.2,-6.2) -- (2.7,-7) -- (3.7,-7);
\draw (3.7,-12.2) node [right] {$baab$};
\draw [black] (21.967,-9.34) arc (142.39859:37.60141:5.217);
\fill [black] (30.23,-9.34) -- (30.14,-8.4) -- (29.35,-9.01);
\draw (26.1,-6.81) node [above] {$abab$};
\draw [black] (30.237,-14.861) arc (-37.53087:-142.46913:5.217);
\fill [black] (21.96,-14.86) -- (22.05,-15.8) -- (22.85,-15.19);
\draw (26.1,-17.4) node [below] {$baba$};
\draw [black] (34.3,-12.1) -- (39.8,-12.1);
\fill [black] (39.8,-12.1) -- (39,-11.6) -- (39,-12.6);
\draw (37.05,-12.6) node [below] {$babb$};
\draw [black] (43.288,-9.154) arc (-198.85998:-280.83693:9.09);
\fill [black] (53.6,-3.16) -- (52.91,-2.52) -- (52.72,-3.5);
\draw (44.89,-3.73) node [above] {$abbb$};
\draw [black] (56.29,-7.2) -- (55.91,-18.2);
\fill [black] (55.91,-18.2) -- (56.43,-17.42) -- (55.43,-17.38);
\draw (55.55,-12.69) node [left] {$bbba$};
\draw [black] (53.178,-22.648) arc (-66.4603:-154.29257:16.023);
\fill [black] (32.34,-14.91) -- (32.24,-15.85) -- (33.14,-15.41);
\draw (38.92,-23.54) node [below] {$bbab$};
\draw [black] (52.861,-21.739) arc (-88.38254:-161.6015:9.786);
\fill [black] (52.86,-21.74) -- (52.05,-21.26) -- (52.08,-22.26);
\draw (44.58,-20.47) node [below] {$abba$};
\end{tikzpicture}
\end{center}

\begin{center}
\begin{tikzpicture}[scale=0.2]
\tikzstyle{every node}+=[inner sep=0pt]
\draw [black] (29.5,-14.5) circle (3);
\draw (29.5,-14.5) node {$abab$};
\draw [black] (17.1,-3.2) circle (3);
\draw (17.1,-3.2) node {$baba$};
\draw [black] (3.2,-13.6) circle (3);
\draw (3.2,-13.6) node {$abaa$};
\draw [black] (9,-28.8) circle (3);
\draw (9,-28.8) node {$baab$};
\draw [black] (22.9,-28.8) circle (3);
\draw (22.9,-28.8) node {$aaba$};
\draw [black] (47.6,-14.5) circle (3);
\draw (47.6,-14.5) node {$babb$};
\draw [black] (60.3,-4) circle (3);
\draw (60.3,-4) node {$abbb$};
\draw [black] (72.9,-14.5) circle (3);
\draw (72.9,-14.5) node {$bbba$};
\draw [black] (61,-14.5) circle (3);
\draw (61,-14.5) node {$abba$};
\draw [black] (60.3,-29.7) circle (3);
\draw (60.3,-29.7) node {$bbab$};
\draw [black] (27.28,-12.48) -- (19.32,-5.22);
\fill [black] (19.32,-5.22) -- (19.57,-6.13) -- (20.25,-5.39);
\draw (26.15,-8.36) node [above] {$ababa$};
\draw [black] (14.7,-5) -- (5.6,-11.8);
\fill [black] (5.6,-11.8) -- (6.54,-11.72) -- (5.94,-10.92);
\draw (7.31,-7.9) node [above] {$babaa$};
\draw [black] (4.27,-16.4) -- (7.93,-26);
\fill [black] (7.93,-26) -- (8.11,-25.07) -- (7.18,-25.43);
\draw (5.35,-22.04) node [left] {$abaab$};
\draw [black] (12,-28.8) -- (19.9,-28.8);
\fill [black] (19.9,-28.8) -- (19.1,-28.3) -- (19.1,-29.3);
\draw (15.95,-29.3) node [below] {$baaba$};
\draw [black] (24.16,-26.08) -- (28.24,-17.22);
\fill [black] (28.24,-17.22) -- (27.45,-17.74) -- (28.36,-18.16);
\draw (26.92,-22.68) node [right] {$aabab$};
\draw [black] (32.5,-14.5) -- (44.6,-14.5);
\fill [black] (44.6,-14.5) -- (43.8,-14) -- (43.8,-15);
\draw (38.55,-15) node [below] {$ababb$};
\draw [black] (48.186,-11.567) arc (161.10807:98.05789:11.319);
\fill [black] (57.31,-4.02) -- (56.45,-3.64) -- (56.59,-4.63);
\draw (48.73,-6.01) node [above] {$babbb$};
\draw [black] (63.249,-4.519) arc (73.86842:26.52044:13.955);
\fill [black] (71.86,-11.69) -- (71.95,-10.75) -- (71.05,-11.2);
\draw (71.2,-6.71) node [above] {$abbba$};
\draw [black] (73.037,-17.49) arc (-3.93328:-75.3806:13.114);
\fill [black] (63.26,-29.28) -- (64.16,-29.56) -- (63.91,-28.59);
\draw (70.6,-26.39) node [right] {$bbbab$};
\draw [black] (57.309,-29.743) arc (-96.55857:-183.68214:11.636);
\fill [black] (47.03,-17.44) -- (46.48,-18.2) -- (47.47,-18.27);
\draw (49.16,-27.08) node [left] {$bbabb$};
\draw [black] (50.6,-14.5) -- (58,-14.5);
\fill [black] (58,-14.5) -- (57.2,-14) -- (57.2,-15);
\draw (54.3,-15) node [below] {$babba$};
\draw [black] (60.86,-17.5) -- (60.44,-26.7);
\fill [black] (60.44,-26.7) -- (60.97,-25.93) -- (59.98,-25.88);
\draw (60.08,-22.08) node [left] {$abbab$};
\end{tikzpicture}
\end{center}

\begin{center}
\begin{tikzpicture}[scale=0.2]
\tikzstyle{every node}+=[inner sep=0pt]
\draw [black] (27.5,-12.8) circle (3);
\draw (27.5,-12.8) node {$aabab$};
\draw [black] (14.8,-3.2) circle (3);
\draw (14.8,-3.2) node {$ababa$};
\draw [black] (3.2,-12.8) circle (3);
\draw (3.2,-12.8) node {$babaa$};
\draw [black] (8.4,-25) circle (3);
\draw (8.4,-25) node {$abaab$};
\draw [black] (21.1,-25) circle (3);
\draw (21.1,-25) node {$baaba$};
\draw [black] (43.3,-12.8) circle (3);
\draw (43.3,-12.8) node {$ababb$};
\draw [black] (43.3,-40.8) circle (3);
\draw (43.3,-40.8) node {$babba$};
\draw [black] (73,-12.8) circle (3);
\draw (73,-12.8) node {$abbba$};
\draw [black] (58.1,-12.8) circle (3);
\draw (58.1,-12.8) node {$babbb$};
\draw [black] (58.1,-28) circle (3);
\draw (58.1,-28) node {$bbabb$};
\draw [black] (73,-28) circle (3);
\draw (73,-28) node {$bbbab$};
\draw [black] (43.3,-28) circle (3);
\draw (43.3,-28) node {$abbab$};
\draw [black] (25.11,-10.99) -- (17.19,-5.01);
\fill [black] (17.19,-5.01) -- (17.53,-5.89) -- (18.13,-5.09);
\draw (23.99,-7.5) node [above] {$aababa$};
\draw [black] (12.49,-5.11) -- (5.51,-10.89);
\fill [black] (5.51,-10.89) -- (6.45,-10.76) -- (5.81,-9.99);
\draw (6.16,-7.51) node [above] {$ababaa$};
\draw [black] (4.38,-15.56) -- (7.22,-22.24);
\fill [black] (7.22,-22.24) -- (7.37,-21.31) -- (6.45,-21.7);
\draw (5.07,-19.85) node [left] {$babaab$};
\draw [black] (11.4,-25) -- (18.1,-25);
\fill [black] (18.1,-25) -- (17.3,-24.5) -- (17.3,-25.5);
\draw (14.75,-25.5) node [below] {$abaaba$};
\draw [black] (22.49,-22.34) -- (26.11,-15.46);
\fill [black] (26.11,-15.46) -- (25.29,-15.93) -- (26.18,-16.4);
\draw (24.98,-20.05) node [right] {$aabab$};
\draw [black] (30.5,-12.8) -- (40.3,-12.8);
\fill [black] (40.3,-12.8) -- (39.5,-12.3) -- (39.5,-13.3);
\draw (35.4,-13.3) node [below] {$aababb$};
\draw [black] (46.3,-12.8) -- (55.1,-12.8);
\fill [black] (55.1,-12.8) -- (54.3,-12.3) -- (54.3,-13.3);
\draw (50.7,-12.3) node [above] {$ababbb$};
\draw [black] (61.1,-12.8) -- (70,-12.8);
\fill [black] (70,-12.8) -- (69.2,-12.3) -- (69.2,-13.3);
\draw (65.55,-13.3) node [below] {$babbba$};
\draw [black] (73,-15.8) -- (73,-25);
\fill [black] (73,-25) -- (73.5,-24.2) -- (72.5,-24.2);
\draw (72.5,-20.4) node [left] {$abbbab$};
\draw [black] (70,-28) -- (61.1,-28);
\fill [black] (61.1,-28) -- (61.9,-28.5) -- (61.9,-27.5);
\draw (65.55,-27.5) node [above] {$bbbabb$};
\draw [black] (58.1,-25) -- (58.1,-15.8);
\fill [black] (58.1,-15.8) -- (57.6,-16.6) -- (58.6,-16.6);
\draw (58.6,-20.4) node [right] {$bbabbb$};
\draw [black] (58.304,-30.984) arc (-3.75714:-94.53211:11.218);
\fill [black] (46.22,-41.43) -- (46.98,-41.99) -- (47.06,-41);
\draw (57.85,-39.22) node [below] {$bbabba$};
\draw [black] (43.3,-37.8) -- (43.3,-31);
\fill [black] (43.3,-31) -- (42.8,-31.8) -- (43.8,-31.8);
\draw (43.8,-34.4) node [right] {$babbab$};
\draw [black] (46.3,-28) -- (55.1,-28);
\fill [black] (55.1,-28) -- (54.3,-27.5) -- (54.3,-28.5);
\draw (50.7,-28.5) node [below] {$abbabb$};
\end{tikzpicture}
\end{center}

\begin{center}
\begin{tikzpicture}[scale=0.2]
\tikzstyle{every node}+=[inner sep=0pt]
\draw [black] (29.2,-12.8) circle (3);
\draw (29.2,-12.8) node {$baabab$};
\draw [black] (16.5,-3.2) circle (3);
\draw (16.5,-3.2) node {$aababa$};
\draw [black] (4.9,-12.8) circle (3);
\draw (4.9,-12.8) node {$ababaa$};
\draw [black] (10.1,-25) circle (3);
\draw (10.1,-25) node {$babaab$};
\draw [black] (22.8,-25) circle (3);
\draw (22.8,-25) node {$abaaba$};
\draw [black] (45,-12.8) circle (3);
\draw (45,-12.8) node {$aababb$};
\draw [black] (45,-40.8) circle (3);
\draw (45,-40.8) node {$bbabba$};
\draw [black] (74.7,-12.8) circle (3);
\draw (74.7,-12.8) node {$babbba$};
\draw [black] (59.8,-12.8) circle (3);
\draw (59.8,-12.8) node {$ababbb$};
\draw [black] (59.8,-28) circle (3);
\draw (59.8,-28) node {$abbabb$};
\draw [black] (74.7,-28) circle (3);
\draw (74.7,-28) node {$abbbab$};
\draw [black] (45,-28) circle (3);
\draw (45,-28) node {$babbab$};
\draw [black] (74.7,-40.8) circle (3);
\draw (74.7,-40.8) node {$bbbabb$};
\draw [black] (26.81,-10.99) -- (18.89,-5.01);
\fill [black] (18.89,-5.01) -- (19.23,-5.89) -- (19.83,-5.09);
\draw (26.63,-7.5) node [above] {$baababa$};
\draw [black] (14.19,-5.11) -- (7.21,-10.89);
\fill [black] (7.21,-10.89) -- (8.15,-10.76) -- (7.51,-9.99);
\draw (6.97,-7.51) node [above] {$aababaa$};
\draw [black] (6.08,-15.56) -- (8.92,-22.24);
\fill [black] (8.92,-22.24) -- (9.07,-21.31) -- (8.15,-21.7);
\draw (6.77,-19.85) node [left] {$ababaab$};
\draw [black] (13.1,-25) -- (19.8,-25);
\fill [black] (19.8,-25) -- (19,-24.5) -- (19,-25.5);
\draw (16.45,-25.5) node [below] {$babaaba$};
\draw [black] (24.19,-22.34) -- (27.81,-15.46);
\fill [black] (27.81,-15.46) -- (26.99,-15.93) -- (27.88,-16.4);
\draw (26.68,-20.05) node [right] {$abaabab$};
\draw [black] (32.2,-12.8) -- (42,-12.8);
\fill [black] (42,-12.8) -- (41.2,-12.3) -- (41.2,-13.3);
\draw (37.1,-13.3) node [below] {$baababb$};
\draw [black] (48,-12.8) -- (56.8,-12.8);
\fill [black] (56.8,-12.8) -- (56,-12.3) -- (56,-13.3);
\draw (52.4,-13.3) node [below] {$aababbb$};
\draw [black] (62.8,-12.8) -- (71.7,-12.8);
\fill [black] (71.7,-12.8) -- (70.9,-12.3) -- (70.9,-13.3);
\draw (67.25,-13.3) node [below] {$ababbba$};
\draw [black] (74.7,-15.8) -- (74.7,-25);
\fill [black] (74.7,-25) -- (75.2,-24.2) -- (74.2,-24.2);
\draw (74.2,-20.4) node [left] {$babbbab$};
\draw [black] (59.8,-25) -- (59.8,-15.8);
\fill [black] (59.8,-15.8) -- (59.3,-16.6) -- (60.3,-16.6);
\draw (59.3,-20.4) node [left] {$abbabbb$};
\draw [black] (45,-37.8) -- (45,-31);
\fill [black] (45,-31) -- (44.5,-31.8) -- (45.5,-31.8);
\draw (45.5,-34.4) node [right] {$bbabbab$};
\draw [black] (48,-28) -- (56.8,-28);
\fill [black] (56.8,-28) -- (56,-27.5) -- (56,-28.5);
\draw (52.4,-28.5) node [below] {$babbabb$};
\draw [black] (74.7,-31) -- (74.7,-37.8);
\fill [black] (74.7,-37.8) -- (75.2,-37) -- (74.2,-37);
\draw (74.2,-34.4) node [left] {$abbbabb$};
\draw [black] (71.7,-40.8) -- (48,-40.8);
\fill [black] (48,-40.8) -- (48.8,-41.3) -- (48.8,-40.3);
\draw (59.85,-40.3) node [above] {$bbbabba$};
\end{tikzpicture}
\end{center}

\begin{center}
\begin{tikzpicture}[scale=0.2]
\tikzstyle{every node}+=[inner sep=0pt]
\draw [black] (22.9,-3.2) circle (3);
\draw (22.9,-3.2) node {$aababbb$};
\draw [black] (5.9,-3.2) circle (3);
\draw (5.9,-3.2) node {$baababb$};
\draw [black] (5.9,-17.4) circle (3);
\draw (5.9,-17.4) node {$abaabab$};
\draw [black] (5.9,-32.2) circle (3);
\draw (5.9,-32.2) node {$babaaba$};
\draw [black] (22.9,-32.2) circle (3);
\draw (22.9,-32.2) node {$ababaab$};
\draw [black] (46.1,-3.2) circle (3);
\draw (46.1,-3.2) node {$ababbba$};
\draw [black] (46.1,-31.2) circle (3);
\draw (46.1,-31.2) node {$babbabb$};
\draw [black] (75.4,-3.2) circle (3);
\draw (75.4,-3.2) node {$abbbabb$};
\draw [black] (61.2,-18.4) circle (3);
\draw (61.2,-18.4) node {$bbabbba$};
\draw [black] (75.4,-18.4) circle (3);
\draw (75.4,-18.4) node {$bbbabba$};
\draw [black] (46.1,-18.4) circle (3);
\draw (46.1,-18.4) node {$abbabbb$};
\draw [black] (75.4,-31.2) circle (3);
\draw (75.4,-31.2) node {$bbabbab$};
\draw [black] (22.9,-17.4) circle (3);
\draw (22.9,-17.4) node {$aababaa$};
\draw [black] (37,-17.4) circle (3);
\draw (37,-17.4) node {$baababa$};
\draw [black] (61.2,-3.2) circle (3);
\draw (61.2,-3.2) node {$babbbab$};
\draw [black] (25.9,-3.2) -- (43.1,-3.2);
\fill [black] (43.1,-3.2) -- (42.3,-2.7) -- (42.3,-3.7);
\draw (34.5,-3.7) node [below] {$aababbba$};
\draw [black] (75.4,-6.2) -- (75.4,-15.4);
\fill [black] (75.4,-15.4) -- (75.9,-14.6) -- (74.9,-14.6);
\draw (74.9,-10.8) node [left] {$abbbabba$};
\draw [black] (46.1,-28.2) -- (46.1,-21.4);
\fill [black] (46.1,-21.4) -- (45.6,-22.2) -- (46.6,-22.2);
\draw (46.6,-24.8) node [right] {$babbabbb$};
\draw [black] (49.1,-18.4) -- (58.2,-18.4);
\fill [black] (58.2,-18.4) -- (57.4,-17.9) -- (57.4,-18.9);
\draw (53.65,-18.9) node [below] {$abbabbba$};
\draw [black] (75.4,-21.4) -- (75.4,-28.2);
\fill [black] (75.4,-28.2) -- (75.9,-27.4) -- (74.9,-27.4);
\draw (74.9,-24.8) node [left] {$bbbabbab$};
\draw [black] (72.4,-31.2) -- (49.1,-31.2);
\fill [black] (49.1,-31.2) -- (49.9,-31.7) -- (49.9,-30.7);
\draw (60.75,-30.7) node [above] {$bbabbabb$};
\draw [black] (34,-17.4) -- (25.9,-17.4);
\fill [black] (25.9,-17.4) -- (26.7,-17.9) -- (26.7,-16.9);
\draw (29.95,-16.9) node [above] {$baababaa$};
\draw [black] (22.9,-20.4) -- (22.9,-29.2);
\fill [black] (22.9,-29.2) -- (23.4,-28.4) -- (22.4,-28.4);
\draw (22.4,-24.8) node [left] {$aababaab$};
\draw [black] (19.9,-32.2) -- (8.9,-32.2);
\fill [black] (8.9,-32.2) -- (9.7,-32.7) -- (9.7,-31.7);
\draw (14.4,-31.7) node [above] {$ababaaba$};
\draw [black] (5.9,-29.2) -- (5.9,-20.4);
\fill [black] (5.9,-20.4) -- (5.4,-21.2) -- (6.4,-21.2);
\draw (5.4,-24.8) node [left] {$babaabab$};
\draw [black] (5.9,-14.4) -- (5.9,-6.2);
\fill [black] (5.9,-6.2) -- (5.4,-7) -- (6.4,-7);
\draw (6.4,-10.3) node [right] {$abaababb$};
\draw [black] (8.9,-3.2) -- (19.9,-3.2);
\fill [black] (19.9,-3.2) -- (19.1,-2.7) -- (19.1,-3.7);
\draw (14.4,-3.7) node [below] {$baababbb$};
\draw [black] (49.1,-3.2) -- (58.2,-3.2);
\fill [black] (58.2,-3.2) -- (57.4,-2.7) -- (57.4,-3.7);
\draw (53.65,-3.7) node [below] {$ababbbab$};
\draw [black] (64.2,-3.2) -- (72.4,-3.2);
\fill [black] (72.4,-3.2) -- (71.6,-2.7) -- (71.6,-3.7);
\draw (68.3,-3.7) node [below] {$babbbabb$};
\draw [black] (61.2,-15.4) -- (61.2,-6.2);
\fill [black] (61.2,-6.2) -- (60.7,-7) -- (61.7,-7);
\draw (60.7,-10.8) node [left] {$bbabbbab$};
\end{tikzpicture}
\end{center}

\end{document}